\documentclass[11pt]{amsart}
\hoffset= - 0.75 in

\usepackage{amsfonts, amsmath, amssymb, amsthm}
\usepackage{latexsym}
\usepackage{color, graphicx, overpic}
\usepackage{subfigure, url}

\newtheorem{thm}{Theorem}

\newtheorem{cor}[thm]{Corollary}

\begin{document}

\author{Michael Joswig}
\thanks{Partially supported by Deutsche
    Forschungsgemeinschaft: SFB~288 ``Differential Geometry and
    Quantum Physics,'' DFG Research Center ``Mathematics for key
    technologies: Modelling, simulation and optimization of real-world
    processes.''}
\address{Institut f\"ur Mathematik, MA6-2, TU Berlin, 10623 Berlin, Germany}
\email{joswig@math.tu-berlin.de}

\title{Computing Invariants of Simplicial Manifolds}

\begin{abstract}
  This is a survey of known algorithms in algebraic topology with a
  focus on finite simplicial complexes and, in particular, simplicial
  manifolds.  Wherever possible an elementary approach is chosen.
  This way the text may also serve as a condensed but very basic
  introduction to the algebraic topology of simplicial manifolds.

  This text will appear as a chapter in the forthcoming book
  ``Triangulated Manifolds with Few Vertices'' by Frank H. Lutz.
\end{abstract}

\maketitle

\newcommand{\NN}{{\mathbb N}}
\newcommand{\QQ}{{\mathbb Q}}
\newcommand{\RR}{{\mathbb R}}
\newcommand{\CC}{{\mathbb C}}
\newcommand{\Sph}{{\mathbb S}}
\newcommand{\ZZ}{{\mathbb Z}}
\newcommand{\cT}{{\mathcal T}}
\newcommand{\rpp}{\RR\textbf{P}^2}
\newcommand{\cpp}{\CC\textbf{P}^2}
\newcommand{\SetOf}[2]{\left\{{#1}\,\mid\,{#2}\right\}}
\newcommand{\size}{\operatorname{size}}
\newcommand{\sign}{\operatorname{sign}}
\newcommand{\abs}[1]{|{#1}|}
\newcommand\classP{\textbf{\rm P}}
\newcommand\classNP{\textbf{\rm NP}}

\definecolor{gold}{rgb}{1.0,0.843,0.0} 


The purpose of this chapter is to survey what is known about algorithms for the computation of algebraic invariants of
topological spaces.  Primarily, we use finite simplicial complexes as our model of topological spaces; for a discussion
of different views see Section~\ref{inv:sec:data}.

On the way we give explicit definitions or constructions of all invariants presented.  Note that we did not try to
phrase all the results in their greatest generality.  Similarly, we focus on invariants for which actual implementations
exist.  The reader is referred to Bredon's monograph~\cite{MR2000b:55001} for the wider perspective.

For a related survey see Vegter~\cite{HDCG:Vegter}.

\section{Homology}

\subsection{Definitions}
Let $\Delta$ be a finite simplicial complex, and let $\Delta^{(k)}$ be the set of its $k$-dimensional faces.  Assume
that the set of vertices is $[n]$.  This way each $k$-dimensional face, or \emph{$k$-face} for short,
$\sigma\in\Delta^{(k)}$ has a unique representation $\{v_0,v_1,\dots,v_k\}$ as an ordered set with $v_i<v_{i+1}$, for
which we simply write $v_0v_1\cdots v_k$.

Fixing a commutative ring $R$ with unit element, which we write as~$1$, allows to define the \emph{boundary map}
$$\partial(v_0v_1\cdots v_k)=\sum_{i=0}^k(-1)^iv_0\cdots v_{i-1}v_{i+1}\cdots v_k.$$
That is to say the
boundary of a $k$-face is a formal linear combinations of $(k-1)$-faces with coefficients in $R$.  The set of all such
linear combinations is denoted by $C_{k-1}(\Delta;R)$, and its elements are called \emph{$(k-1)$-chains}. The
$k$-dimensional simplicial \emph{boundary operator} with coefficients in~$R$ is the unique $R$-linear extension of the
boundary map to the $R$-module of $k$-chains $C_k(\Delta;R)$; it is denoted by $\partial_k$.

The image of $\partial_{k+1}$, denoted by $B_k(\Delta;R)$, and the kernel of $\partial_k$, denoted by $Z_k(\Delta;R)$,
are both $R$-modules.  Their elements are called \emph{$k$-boundaries} and \emph{$k$-cycles}, respectively.  A direct
computation yields that $\partial_k(\partial_{k+1}(\sigma))=0$ or, equivalently, $B_k(\Delta;R)\le Z_k(\Delta;R)$.  The
quotient $$H_k(\Delta;R)=Z_k(\Delta;R)/B_k(\Delta;R)$$
is the $k$-th \emph{homology module} of~$\Delta$ with
coefficients in~$R$.

There is a certain freedom of choice in defining the map $\partial_0$.  Usually, we take $\partial_0=0$.  Hence
$H_0(\Delta;R)\cong R$ if and only if $\Delta$ is connected.  We call a connected complex \emph{$R$-acyclic} if all
homology modules $H_k(\Delta;R)$ vanish for $k>0$.  Contractible complexes are $R$-acyclic for arbitrary~$R$.  If,
however, we consider the empty set as a simplex of dimension $-1$, then the space of $(-1)$-chains is isomorphic to~$R$
as an $R$-module and $\partial_0$ is surjective.  By definition $\partial_{-1}=0$.  This leads to what is called
\emph{reduced homology}: $\widetilde H_k(\Delta;R)=H_k(\Delta;R)$ for $k>0$ and $H_0(\Delta;R)\cong\widetilde
H_0(\Delta;R)\oplus R$.

The \emph{dimension} of $\Delta$, denoted as $\dim\Delta$, is the maximal dimension of its faces.  Clearly, for
$k>\dim\Delta$, we have $C_k(\Delta;R)=0$ and hence also $H_k(\Delta;R)=0$.

There are several other homology theories.  Among the most important are singular and cellular homology.  All these
theories coincide for ``reasonable'' spaces including CW-complexes and, in particular, simplicial complexes;
see~\cite[Chapter~IV]{MR2000b:55001}.

Not all coefficient domains are equally important. The Universal Coefficient Theorem~\cite[V.7.5]{MR2000b:55001} implies
that knowing the integral homology of a space, that is, the case where $R=\ZZ$, determines the homology with many other
coefficient domains, including the rationals and all finite fields.  Conversely, the integral homology is determined by
the rational homology together with the homology over all finite fields of prime order.  Usually, we write $H_k(\Delta)$
for $H_k(\Delta;\ZZ)$ and likewise for other $\ZZ$-modules.

\subsection{Algorithms and Complexity}

It is fairly easy to actually compute the homology of a space if the coefficient domain is a field.  In this case the
boundary operator is a linear map.  Its kernel and image are vector spaces.  Hence the homology modules are vector
spaces whose dimensions can be obtained as the differences of co-rank and rank, respectively, of two subsequent boundary
matrices.  Algorithmically, the rank (or co-rank) of a matrix can be obtained by a sequence of Gauss elimination steps.

We denote the number of $k$-faces of~$\Delta$ by $f_k(\Delta)$ or simply $f_k$.  It follows that for a field~$F$ it
requires $O(\max(f_{k+1},f_k,f_{k-1})^3)$ arithmetic operations in~$F$ to compute $H_k(\Delta;F)$ \emph{provided that}
the boundary matrices $\partial_{k+1}$ and~$\partial_k$ are given.  This latter condition can be read as a warning:
Typically a finite simplicial complex is given as the set of faces which are maximal with respect to inclusion.  Even if
their number is small (e.g., one for a simplex) the lower dimensional boundary matrices can still be huge: The trivial
upper bound $f_k\le\binom{n}{k}$ is tight.

We define the \emph{size} of $\Delta$ to be the number of incidences between vertices and facets.  In particular, the
size is bounded from above by $(\dim\Delta+1)m$, where $m$ is the number of facets of~$\Delta$.  Equality holds if and
only if $\Delta$ is pure.  This notion corresponds to an encoding of $\Delta$ as an array of its facets, where each
facet is represented as a list of its vertices.

Each integral homology module of~$\Delta$ is a finitely generated abelian group, that is to say, a direct product of a
finitely generated free abelian group and a product of finitely many cyclic groups of prime power order.  The subgroup
of elements of finite order is called the \emph{torsion} subgroup.  We will see below that every finitely generated
abelian group, in fact, arises as a homology module of some $2$-dimensional simplicial complex.

On a conceptual level the computation of integral homology differs only by a little from the computation of rational
homology.  There are two reasons for this: Firstly, $\QQ$ is the quotient field of~$\ZZ$.  Secondly, $\ZZ$ is a
Euclidean domain.  The latter property comes in handy when one tries to perform a Gauss elimination over~$\ZZ$. Either
the pivoting element is a unit, that is $1$ or $-1$, then the corresponding rational basis transformations are integral,
and the situation is exactly the same as over~$\QQ$.  Or, if the pivoting element is not a unit, Euclid's algorithm can
be used to determine an integral transformation which simplifies the matrix.  In this way, after finitely many steps,
each boundary matrix $\partial_k$ can be transformed into diagonal shape, the \emph{Smith-Normal-Form} of~$\partial_k$.
It is then straightforward to determine the structure of $H_k(\Delta)$ from the Smith-Normal-Forms of $\partial_{k+1}$
and $\partial_k$.

The existence and also the construction of the Smith-Normal-Form of an integral matrix is a classical result due to H.
J. S. Smith~\cite{SNF}; for a modern account see Munkres~\cite[\S11]{MR85m:55001}.

In terms of complexity there is a fundamental difference between Gauss elimination over~$\QQ$ and the sketched type of
Euclidean-Gaussian elimination over~$\ZZ$.  In the rational case the computation time is bounded by a polynomial (in the
size of the input matrix) since the growth of the coefficients can be controlled.  This has been observed by
Edmonds~\cite{MR37:5114}; see Schrijver~\cite[\S3.3]{MR88m:90090}.  However, for the analogous operation over~$\ZZ$
neither the size of the coefficients nor the number of arithmetic operations in~$\ZZ$ is polynomially bounded.

Relying on modular arithmetic, though, Kannan and Bachem~\cite{MR81k:15002} gave a first polynomial time
Smith-Normal-Form algorithm, which was later improved by Iliopoulos~\cite{MR91a:20065}, Storjohann~\cite{MR99j:65074},
and others.

While the modular approach is valid for matrices with arbitrary integer coefficients, simplicial boundary matrices have
entries $1$, $-1$, and $0$ only.  That is to say, in an arbitrary simplicial boundary matrix it is \emph{always}
possible to perform at least a few Gauss elimination steps.  Moreover, a typical boundary matrix is very sparse.  If the
matrix stays sparse during the elimination and if, additionally, one does not run out of unit coefficients too early
(such that it is possible to continue with elimination steps) an elimination based Smith-Normal-Form algorithm can be
superior to the more sophisticated methods.  This is why in practical applications elimination algorithms are often
preferred; see the survey of Dumas, Heckenbach, Saunders, and Welker~\cite{AGSS:DumasEtAl} and also the
Section~\ref{inv:subsec:morse} on discrete Morse theory.

Donald and Chang~\cite{MR93j:68086} proved that, asymptotically, for random simplicial complexes with many vertices
compared to the dimension~$d$ the integral homology can be computed in $O(df_{\textrm{max}}^2)$ time with high
probability.  Here $f_{\textrm{max}}=\max\{f_0,\dots,f_d\}$.

\subsection{Examples}\label{inv:homology:examples}

\begin{figure}
  \centering
  \begin{overpic}[width=5.5cm]{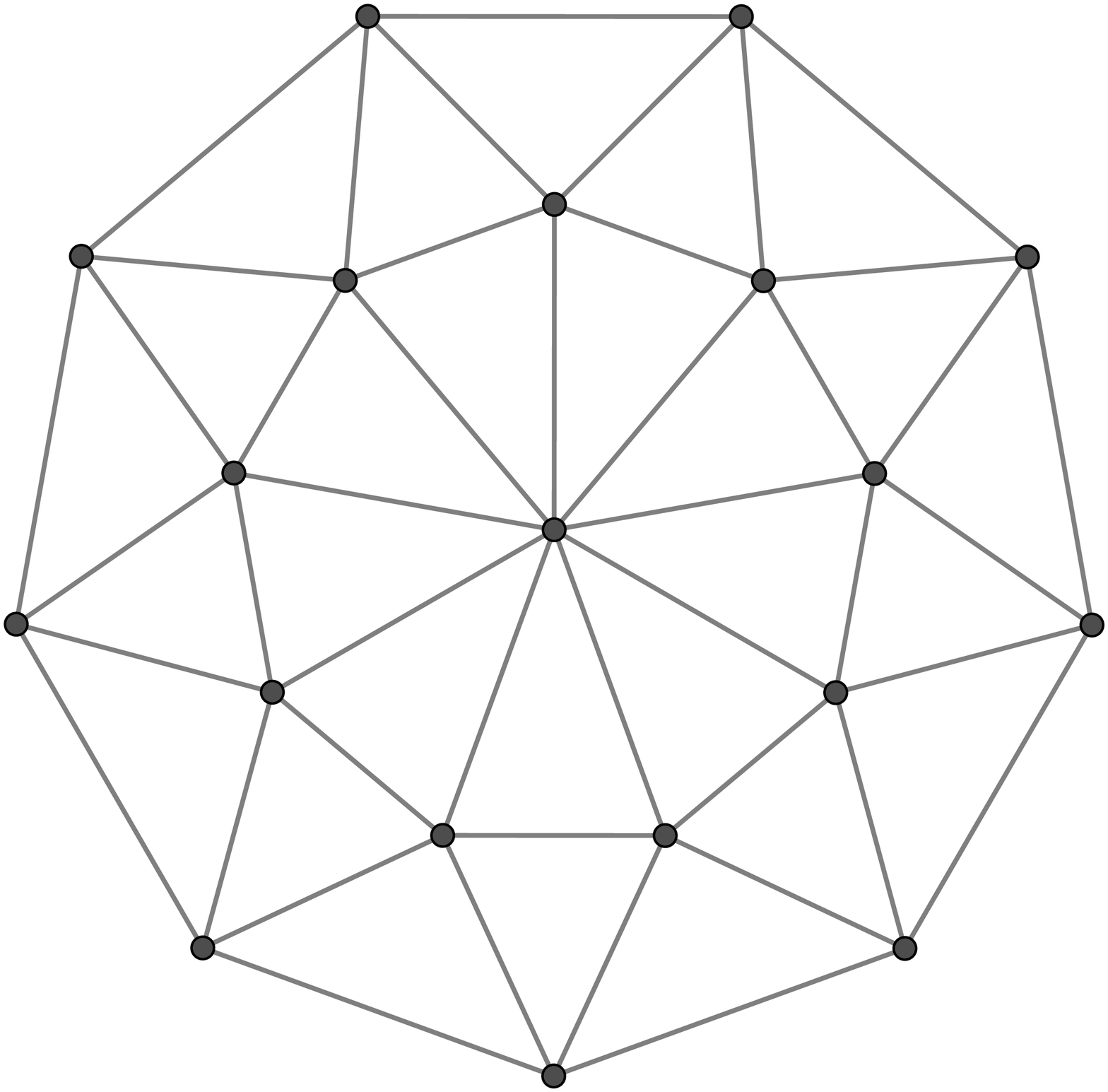}
    \put(48,-3){$a$}
    \put(80,10){$b$}
    \put(96,40){$c$}
    \put(90,73){$a$}
    \put(65,94){$b$}
    \put(30,94){$c$}
    \put(5,73){$a$}
    \put(-1,40){$b$}
    \put(15,10){$c$}
  \end{overpic}
  \caption{Complex $C(3)$: All unlabeled vertices are distinct.  The three special vertices $a,b,c$ occur three times
    each in the boundary of the $2$-cell suggested by the picture.\label{inv:fig:c3}}
\end{figure}

The homology of a $1$-dimensional simplicial complex, that is, a graph, is very easy to understand. For a connected
graph $\Gamma$ with $n$ vertices and $e$ edges we immediately obtain $H_1(\Gamma)\cong\ZZ^{e-n+1}$: Consider a spanning
tree, which has $n-1$ edges; each one of the $e-(n-1)$ non-tree edges then corresponds to a different free generator.
However, for each finitely generated abelian group $G$ there is a $2$-dimensional simplicial complex~$\Delta$ such that
$H_1(\Delta)\cong G$: There is a family of connected $2$-dimensional simplicial complexes $C(k)$, for $k\ge 2$, with
$f(C(k))=(1+3(k+1),3+12k,9k)$ such that $H_2(C(k))=0$ and $H_1(C(k))\cong\ZZ/k$.  From Figure~\ref{inv:fig:c3}, which
displays~$C(3)$, it should be easy to deduce how to define $C(k)$ for arbitrary~$k$; the complex $C(2)$ is homeomorphic
to the real projective plane~$\rpp$ and is the only manifold in the series.  The non-trivial homology cycle is carried
by the induced subcomplex with vertices $\{a,b,c\}$.  Clearly, $C(k)$ is $\QQ$-acyclic.  A \emph{wedge}
$\Delta\vee\Delta'$ of two simplicial complexes $\Delta,\Delta'$ is a complex arising by identifying one vertex
in~$\Delta$ with one vertex in~$\Delta'$ in the disjoint union $\Delta\sqcup\Delta'$.  It is immediate that
$H_k(\Delta\vee\Delta';R)\cong H_k(\Delta;R)\oplus H_k(\Delta';R)$ for $k\ge1$.  By successively wedging complexes
$C(k)$ and (triangulated) cylinders~$\Sph^1\times[0,1]$ we obtain each finitely generated abelian group as the first
integral homology group of a pure $2$-dimensional simplicial complex.

\subsection{Elimination Strategies and Discrete Morse Theory}\label{inv:subsec:morse}

Classical Morse theory allows to bound the size of the homology modules of a smooth manifold by counting critical points
(with respect to some smooth function) in a cell decomposition; see Milnor~\cite[\S{}I.5]{MR29:634}.
Forman~\cite{MR99b:57050,MR2000h:57041} developed a Morse theory for abstract simplicial complexes (and, more generally,
for arbitrary finite CW-complexes), which we briefly sketch now.

Let $\Delta$ be a finite simplicial complex, and let $\mu:\Delta\to\NN$ be an arbitrary function.  For each simplex
$\sigma\in\Delta$ we define the \emph{critical defect} $$\delta_\mu(\sigma)=\#\SetOf{\rho\in\Delta}{\rho\lessdot\sigma,\ 
  \mu(\rho)\ge\mu(\sigma)}+\#\SetOf{\tau\in\Delta}{\tau\gtrdot\sigma,\ \mu(\tau)\le\mu(\sigma)},$$
where by $\rho\lessdot\sigma$ we denote that $\rho$ is a maximal face of $\sigma$.  Now a function $\mu:\Delta\to\NN$ is
a \emph{discrete Morse function} of~$\Delta$ if for each simplex the critical defect is at most one.  A simplex of
critical defect zero is called \emph{critical}; otherwise it is \emph{regular}.  Discrete Morse functions always exist:
The function which maps each simplex to its dimension is a Morse function for which each simplex is critical.

\begin{figure}[htbp]
  \centering
  \includegraphics[width=7cm]{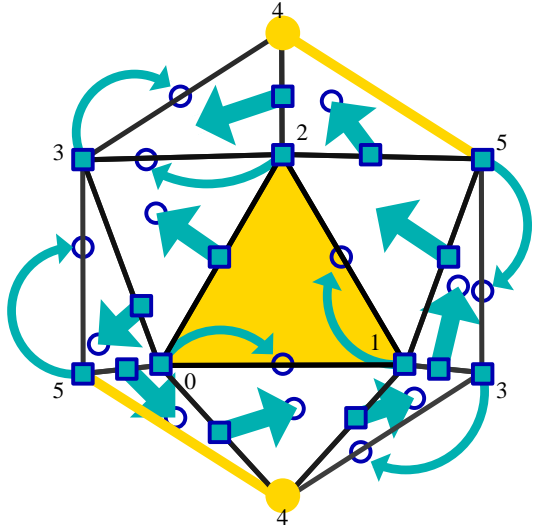}
  \caption{Triangulation $\rpp_6$ of the projective plane (obtained from identifying antipodal points in the
    boundary of the regular icosahedron) with a Morse matching~$[\mu]$ (indicated by the arrows).  This triangulation is
    minimal with $f$-vector $(6,15,10)$.  The regular faces are marked with a circle or a square, depending whether they
    are matched to higher or lower dimensional face.  There are three critical faces: $4$, $45$, and $012$.  It is clear
    that $\mu$ cannot be extended: The critical triangle $012$ cannot be matched since all its edges are already
    matched; the critical vertex~$4$ cannot be matched to the critical edge~$45$ since this would violate the acyclicity
    condition. The subsequent Figures \ref{inv:fig:rpp:matching} and \ref{inv:fig:rpp:elimination} display the same
    information in different ways.\label{inv:fig:rpp:flow}}
\end{figure}

\begin{figure}[htbp]
  \centering
  \includegraphics[width=11cm]{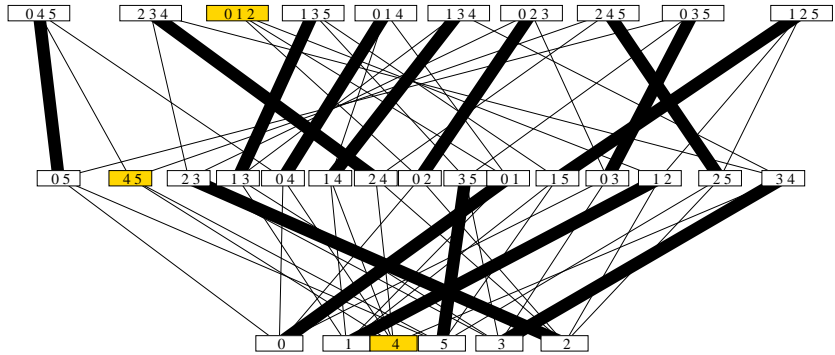}
  \caption{The Morse matching~$[\mu]$ in the Hasse diagram of the face lattice of
    $\rpp_6$.  The bold edges correspond to the arrows in Figure~\ref{inv:fig:rpp:flow}.  The critical faces are
    unmatched.\label{inv:fig:rpp:matching}}
\end{figure}

Fix a discrete Morse function $\mu$.  Then, by definition, for each regular simplex $\sigma$ there is either a unique
regular face $\rho\lessdot\sigma$ with $\mu(\rho)\ge\mu(\sigma)$ or a unique regular face $\tau\gtrdot\sigma$ with
$\mu(\tau)\le\mu(\sigma)$.  That is to say, each discrete Morse function gives rise to a matching in the Hasse diagram
of the face lattice of~$\Delta$, the \emph{Morse matching} of~$\mu$.  Each Morse matching $[\mu]$ satisfies the
following property: Direct all edges of the Hasse diagram of~$\Delta$ consistently, say, with increasing dimension, and
reverse the orientation for each edge in~$[\mu]$; then there is no oriented cycle.  Conversely, each matching satisfying
this acyclicity condition is associated with a discrete Morse function.  For an example of a Morse matching in the
triangulation~$\rpp_6$ of the real projective plane see Figures \ref{inv:fig:rpp:flow}, \ref{inv:fig:rpp:matching}, and
\ref{inv:fig:rpp:elimination}.

Each edge $\sigma\lessdot\tau$ in the Hasse diagram of the face lattice of~$\Delta$, for $\tau\in\Delta^{(k)}$ (and thus
$\sigma\in\Delta^{(k-1)}$), corresponds to a $\pm1$ entry in the boundary matrix $\partial_k$ and hence can be used to
eliminate all other entries in the row of~$\tau$ and in the column of~$\sigma$.  Now, a subset $E$ of the edges of the
Hasse diagram corresponds to a sequence of elimination steps which can be performed (in an arbitrary order) without
interfering with each other if and only if $E$ is a Morse matching.  That is to say, a globally \emph{optimal}
elimination strategy corresponds to a discrete Morse function with as few critical faces as possible.

\begin{figure}[htbp]
  $$
  \begin{array}{c|ccccccccccccccc}
    \partial_2 & 01 & 02 & 03 & 04 & 05 & 12 & 13 & 14 & 15 & 23 & 24 & 25 & 34 & 35 & \colorbox{gold}{45} \\
    \hline
    \colorbox{gold}{012}        &  1 & -1 &    &    &    &  1 &    &    &    &    &    &    &    &    &    \\
    014        &  1 &    &    & \framebox{-1} &    &    &    &  1 &    &    &    &    &    &    &    \\
    023        &    &  \framebox{1} & -1 &    &    &    &    &    &    &  1 &    &    &    &    &    \\
    035        &    &    &  \framebox{1} &    & -1 &    &    &    &    &    &    &    &    &  1 &    \\
    045        &    &    &    &  1 & \framebox{-1} &    &    &    &    &    &    &    &    &    &  1 \\
    125        &    &    &    &    &    &  1 &    &    & \framebox{-1} &    &    &  1 &    &    &    \\
    134        &    &    &    &    &    &    &  1 & \framebox{-1} &    &    &    &    &  1 &    &    \\
    135        &    &    &    &    &    &    &  \framebox{1} &    & -1 &    &    &    &    &  1 &    \\
    234        &    &    &    &    &    &    &    &    &    &  1 & \framebox{-1} &    &  1 &    &    \\
    245        &    &    &    &    &    &    &    &    &    &    &  1 & \framebox{-1} &    &    &  1 

    \\
    \\

    (\partial_1)^{\text{tr}}
               & 01 & 02 & 03 & 04 & 05 & 12 & 13 & 14 & 15 & 23 & 24 & 25 & 34 & 35 & \colorbox{gold}{45} \\
    \hline
    0          & \framebox{-1} & -1 & -1 & -1 & -1 &    &    &    &    &    &    &    &    &    &    \\    
    1          &  1 &    &    &    &    & \framebox{-1} & -1 & -1 & -1 &    &    &    &    &    &    \\ 
    2          &    &  1 &    &    &    &  1 &    &    &    & \framebox{-1} & -1 & -1 &    &    &    \\ 
    3          &    &    &  1 &    &    &    &  1 &    &    & -1 &    &    & \framebox{-1} & -1 &    \\ 
    \colorbox{gold}{4}          &    &    &    &  1 &    &    &    &  1 &    &    &  1 &    & -1 &    &  1 \\ 
    5          &    &    &    &    &  1 &    &    &    &  1 &    &    &  1 &    &  \framebox{1} &  1 
  \end{array}
  $$
  \caption{Boundary matrices $\partial_2$ and $\partial_1$ of $\rpp_6$.  The matrix $\partial_1$ is  transposed
    in order to visualize how the two matrices fit together.  The boxed entries in the two matrices correspond to the
    arrows in Figure~\ref{inv:fig:rpp:flow} and the matching edges in Figure~\ref{inv:fig:rpp:matching}.  Using the
    boxed entries as pivots in an arbitrary order yields an elimination strategy for the computation of the homology
    modules.\label{inv:fig:rpp:elimination}}
\end{figure}

Figure~\ref{inv:fig:rpp:elimination} shows the Morse matching~$[\mu]$ of $\rpp_6$ from Figures \ref{inv:fig:rpp:flow} and
\ref{inv:fig:rpp:matching} as a set of entries in the boundary matrices.  Interpreting each elimination step
geometrically leads to a transformation from a triangulation to a more general finite CW-complex.  In this particular
case, we obtain a cell decomposition of $\rpp$ into three cells or dimensions $0,1,2$, respectively.  The unique one
cell corresponds to a (real) projective line.  This cell decomposition is minimal since a CW-complex with only two cells
(one of which must be a point, by definition) necessarily is a sphere.  Further this argument implies that the Morse
matching~$\mu$ of $\rpp_6$ (which leads to this minimal cell decomposition) is an optimal Morse matching.  Note that
Forman more generally proved \cite[Corollary 1.7]{MR2000h:57041} that a combinatorial manifold without boundary which
admits a Morse matching with only two critical cells is a shellable sphere.  In this sense optimal Morse matchings
generalize shelling orders to arbitrary simplicial complexes.

While optimal discrete Morse functions are important, it is not easy to find good discrete Morse functions
algorithmically.

\begin{thm}{\rm (E\u{g}ecio\u{g}lu and Gonzalez \cite{MR97a:68159})}\label{inv:thm:erasable}
  For the class of finite pure $2$-dimensional simplicial complexes which are embeddable into~$\RR^3$ it is
  \classNP-complete to decide whether there is a Morse matching with at most $c$ critical $2$-faces.
\end{thm}

In fact, E\u{g}ecio\u{g}lu and Gonzalez even proved a rather strong non-approxi\-ma\-bi\-li\-ty result, which we omit.
Notice that, for the result above and the corollary below, the number~$c$ is part of the input.

One can show that an optimal discrete Morse function of any connected simplicial complex necessarily has a unique
critical vertex.  This can be used to obtain the following result~\cite{Complexity+of+Discrete+Morse+Theory}; see also
Lewiner~\cite{Lewiner:Thesis}.

\begin{cor} For the class of finite pure $2$-dimensional simplicial complexes which are embeddable into~$\RR^3$ it is
  \classNP-complete to decide whether there is a Morse matching with at most $c$ critical faces of arbitrary dimension.
\end{cor}

For a more stringent interpretation of discrete Morse theory from the viewpoint of shellings see
Chari~\cite{MR2001g:52016}.  The set of all optimal Morse functions of a given complex may have an interesting
structure, although in most cases it is very difficult to obtain; see \cite{DMC} for a discussion of special cases.

In practice a greedy choice of an elimination strategy/discrete Morse function often leads to a considerable reduction
of the input to a homology computation which is then passed on to some Smith-Normal-Form algorithm.  This way it is
possible to compute the homology of rather large simplicial complexes; e.g., see~\cite{AGSS:DumasEtAl}.

\subsection{Bases Transformations}

In order to locate subcomplexes which carry the homology cycles it is necessary to keep track of all the transformations
of bases during an elimination process.  Even if the boundary matrix stays sparse, the accumulated bases transformation
matrices almost never do.  If we want to avoid to run out of memory during the computation it is therefore useful to
store the transformations as sequences of elementary transformations, which in turn are stored in compressed form.  Such
techniques are also common in numerical linear algebra.

\section{Cohomology}

Transposing each boundary matrix leads to coboundary maps, cochains, cocycles, and, finally, cohomology modules, which
contain the same information as the homology modules, but in a dual form.  Everything concerning the computation of
homology modules can directly be translated to cohomology.  Here we are interested in the cup product, which has a less
natural description in terms of homology.

\subsection{Definitions}

Let again $\Delta$ be a finite simplicial complex, and let $R$ be a commutative ring.  The $R$-module of $k$-dimensional
\emph{cochains} is the set $C^k(\Delta;R)$ of $R$-linear maps from $C_k(\Delta)$ to~$\RR$.  The evaluation of $f\in
C^k(\Delta;R)$ at $\sigma\in C_k(\Delta;R)$ is written as $\langle f,\sigma\rangle$.  For $\sigma=v_0\cdots
v_k\in\Delta^{(k)}$ let $\sigma^*$ be the $k$-cochain with maps $\sigma$ to~$1$ and all other $k$-simplices to~$0$, and
for each face $\tau=v_0\cdots v_i v v_{i+1}\cdots v_k\in\Delta^{(k+1)}$ let $\epsilon(\sigma,\tau)=(-1)^{i+1}$.  The
latter value is the coefficient of $\sigma$ in the boundary of~$\tau$.  We define the \emph{coboundary} of~$\sigma^*$ as
$$\delta\sigma^*=\sum_{\tau\gtrdot\sigma}\epsilon(\sigma,\tau)\tau^*.$$
Via the canonical identification
$\sigma\leftrightarrow\sigma^*$ we have that the $R$-linear extension $\delta^k:C^k(\Delta;R)\to C^{k+1}(\Delta;R)$ of
$\delta$ is the transpose of the map~$\partial_{k+1}$.

Again, the image of $\delta^{k-1}$, denoted by $B^k(\Delta;R)$, and the kernel of $\delta^k$, denoted by
$Z^k(\Delta;R)$, are both $R$-modules.  Their elements are called \emph{$k$-coboundaries} and \emph{$k$-cocycles},
respectively.  The quotient $$H^k(\Delta;R)=Z^k(\Delta;R)/B^k(\Delta;R)$$
is the $k$-th \emph{cohomology module}
of~$\Delta$ with coefficients in~$R$.

One reason to look at cohomology rather than homology is that cohomology comes equipped with a natural $R$-bilinear
function, the \emph{cup product} $\cup:H^i(\Delta;R)\times H^k(\Delta;R)\to H^{i+k}(\Delta;R)$ which is defined on the
cochain level by
$$\langle f\cup g,v_0\cdots v_{i+k}\rangle=(-1)^{ik}\,\langle f,v_0\cdots v_i\rangle\,\langle
g, v_i\cdots v_{i+k}\rangle,$$
where $f\in C^i(\Delta;R)$ and $g\in C^k(\Delta;R)$.  This turns
$H^*(\Delta;R)=\bigoplus_k H^k(\Delta;R)$ into an associative $R$-algebra, the \emph{cohomology algebra}
of~$\Delta$.  Note that the cup product is not commutative, in general. Instead we have
\begin{equation}\label{inv:eq:cup}
g\cup f=(-1)^{ik}f\cup g.
\end{equation}

There is another interesting product, the \emph{cap product} $\cap:H^i(\Delta;R)\times H_k(\Delta;R)\to
H_{k-i}(\Delta;R)$, which can be defined, again on the level of cochains and chains, by
$$f\cap v_0\cdots v_k=(-1)^{i(k-i)}\,\langle f, v_0\cdots v_i\rangle\, v_i\cdots v_k,$$
for $i\le k$.  For details
see~\cite[\S{}VI.5]{MR2000b:55001}.

\subsection{Algorithms}

For the integral homology modules it was obvious what the output of a computation should be: Each finitely generated
abelian group has a canonical finite representation by its free rank and its torsion coefficients.  The goal for
computing the cup product is somewhat less obvious: Computing the cohomology modules is the same as for homology.  If we
keep track of all basis transformations we also obtain a set of cochains which generates $H^*(\Delta;R)$ as an
$R$-module.  This in turn can be used to compute a multiplication table for the cup product; see
Peitgen~\cite{MR51:4208}.

However, there is no canonical choice for a generating system.  This has the serious disadvantage that there does not
seem to be a direct way to decide whether two cohomology rings are isomorphic or not.

\subsection{Duality of Manifolds}

Let $M$ be a closed $d$-dimensional simplicial manifold.  Since each co\-di\-men\-sion-$1$-face is contained in
exactly two facets we have $H_d(M;\ZZ/2)\cong\ZZ/2$.  We denote the generator of $H_d(M;\ZZ/2)$ by $[M]_2$.  The
manifold $M$ is \emph{oriented} if $H_d(M)\cong\ZZ$.  If $M$ is not oriented then $H_d(M)=0$.

The following theorem comes in many guises:

\begin{thm}{\rm (Poincar\'e Duality; see~\cite[VI.8.4]{MR2000b:55001})}\label{inv:thm:Poincare}
  The map $$H^k(M;\ZZ/2)\to H_{d-k}(M;\ZZ/2):f\mapsto f\cap[M]_2$$
  is an isomorphism.  If, additionally, $M$ is
  oriented, and if $[M]$ is a generator of $H_d(M)$, then the map $$H^k(M)\to H_{d-k}(M):f\mapsto f\cap[M]$$
  is an
  isomorphism.
\end{thm}

\subsection{Stiefel-Whitney Classes}

Let $M$ be a closed $d$-dimensional smooth combinatorial manifold. The $k$-th Stiefel-Whitney class of $M$ is
a certain (characteristic) cohomology class $\omega^k\in H^k(M;\ZZ/2)$.  For a precise definition
see~\cite[\S{}VI.17]{MR2000b:55001}.  The $k$-th Stiefel-Whitney class forms the essential obstruction to the existence
of $d-k+1$ linearly independent cross sections in the tangent bundle of~$M$; see Milnor and
Stasheff~\cite[Proposition 4.4]{MR55:13428}.  In the sequel we give an elementary combinatorial description of
$\omega^k$ due to Goldstein and Turner~\cite{MR54:3724}; see also Halperin and Toledo~\cite{MR47:1072}.

For a pair of faces $\sigma=v_0\cdots v_k$ and $\tau$ with $\sigma\leq\tau$ define the $i$-th \emph{gap} as
$$G_i(\sigma,\tau)=\SetOf{w\in\tau}{v_i<w<v_{i+1}},$$ where $0\le i<k$.  Additionally, let
$G_{-1}(\sigma,\tau)=\SetOf{w\in\tau}{w<v_0}$ and $G_k=\SetOf{w\in\tau}{v_k<w}$.  Now $(\sigma,\tau)$ is a \emph{regular
  pair} if for all $i$ odd $G_i(\sigma,\tau)=\emptyset$.

Let $\omega_k$ be the $k$-dimensional chain with coefficients in~$\ZZ/2$ which is the sum of all $k$-faces $\sigma$, for
which the number of $l$-faces $\tau$, such that $l\ge k$ and $(\sigma,\tau)$ is a regular pair, is odd.  It is easy to
verify that $\omega_k$ is a homology cycle mod~$2$.  We identify $\omega_k$ with its $\ZZ/2$-homology class.

Since for each $k$-face $\sigma$ the pair $(\sigma,\sigma)$ is regular, it is obvious that $\omega_d=[M]_2\ne 0$.

\begin{thm}{\rm (Goldstein and Turner~\cite{MR54:3724})}
  The $\ZZ/2$-homology class of $\omega_k$ is the image of the $(d-k)$-th Stiefel-Whitney class $\omega^{d-k}$ under the
  Poincar\'e Duality isomorphism.
\end{thm}

\subsection{The Intersection Form of a $4$-Manifold}

Let $M$ be a closed oriented combinatorial $4$-manifold.  This implies that $M$ carries a unique smooth
structure.

Since $M$ is oriented we have $H_4(M)\cong\ZZ$.  We fix a generator $[M]\in H_4(M)$, which is called
an \emph{orientation class}.  In view of equation~(\ref{inv:eq:cup}) the map
$$Q_M:H^2(M)\times H^2(M)\to \ZZ:(f,g)\mapsto \langle f\cup g,[M]\rangle$$
is a symmetric bilinear form, the
\emph{intersection form} of~$M$.  If $f$ or $g$ is of finite order, then $Q_M(f,g)=0$.  Therefore, by choosing a basis,
we can identify $Q_M$ with an integral $(r\times r)$-matrix $M$, the \emph{Gram matrix} of~$Q_M$, where $r$ is the free
rank of $H^2(M)$.  It can be shown that $Q_M$ is unimodular, that is, $\det Q_M=\pm 1$.  The form $Q_M$ is \emph{even}
if all diagonal elements of~$M$ are even and \emph{odd} otherwise.  Note that this property does not depend on the basis
chosen.  Over the reals $Q_M$ can be transformed into a diagonal matrix $D$ with entries $\pm1$.  The number
$\sign(Q_M)=(\text{number of $1$-entries})-(\text{number of $-1$-entries})\text{ in $D$}$ is the \emph{signature}
of~$Q_M$.  If $Q_M$ is indefinite, that is, $\abs{\sign(Q_M)}\ne r$, then up to an integral basis transformation $Q_M$
is uniquely determined by the rank, the parity, and the signature; see Milnor and Husemoller~\cite[Theorem
II.5.3]{MR58:22129}.  The situation is more complicated if $Q_M$ is a (positive or negative) definite form; see
\cite[\S{}II.6]{MR58:22129}.

The Wu formula~\cite[Proposition 1.4.18]{MR2000h:57038} implies that the intersection form of $M$ is even if and
only if the second Stiefel-Whitney class vanishes.

Whitehead proved that the homotopy type of a simply connected closed $4$-manifold is determined by its intersection
form; see \cite[Theorem V.1.5]{MR58:22129}.  However, this has been improved to a considerably stronger famous result.

\begin{thm}{\rm (Freedman~\cite{MR84b:57006})}\label{inv:thm:Freedman}
  Two simply connected closed oriented combinatorial $4$-manifolds are homeomorphic if and only if their intersection
  forms are equivalent.
\end{thm}

If we additionally assume that $M$ is simply connected, then $H_1(M)=0$ and hence, by the Universal Coefficient
Theorem~\cite[V.7.5]{MR2000b:55001}, $H^1(M)=0$.  The Poincar\'e Duality Theorem~\ref{inv:thm:Poincare} now gives that
$H_3(M)=0$.  Moreover, $H_2(M)\cong H^2(M)$ is torsion-free and thus the intersection form $Q_M$ completely determines
the cohomology ring of~$M$.

The intersection form can also be defined in terms of the cap product (see~\cite[Corollary~VI.5.3]{MR2000b:55001})
\begin{equation}
  \langle f\cup g,[M]\rangle = \langle f,g\cap[M]\rangle,
\end{equation}
where $f,g\in H^2(M)$.  Recall that the map $g\mapsto g\cap[M]$ is the Poincar\'e Duality isomorphism.

For more information on the topology of $4$-manifolds see Kirby~\cite{MR90j:57012}, Gompf and Stipsicz~\cite{MR2000h:57038}.

\subsection{An Example}\label{inv:subsec:cpp}

The following is the list of facets of the unique minimal triangulation $\cpp_9$ of the complex projective plane with
$f$-vector $(9,36,84,90,36)$; see K\"uhnel and Banchoff~\cite{MR85d:51016}:
{\arraycolsep=.15cm%
$$
  \begin{array}{ccccccccc}
    {0 1 2 3 4}
    &{0 1 2 3 6}
    &{0 1 2 4 7}
    &{0 1 2 6 7}
    &{0 1 3 4 5}
    &{0 1 3 5 6}
    &{0 1 4 5 7}
    &{0 1 5 6 8}
    &{0 1 5 7 8}\\
    {0 1 6 7 8}
    &{0 2 3 4 8}
    &{0 2 3 6 7}
    &{0 2 3 7 8}
    &{0 2 4 5 7}
    &{0 2 4 5 8}
    &{0 2 5 7 8}
    &{0 3 4 5 6}
    &{0 3 4 6 8}\\
    {0 3 6 7 8}
    &{0 4 5 6 8}
    &{1 2 3 4 8}
    &{1 2 3 5 6}
    &{1 2 3 5 8}
    &{1 2 4 6 7}
    &{1 2 4 6 8}
    &{1 2 5 6 8}
    &{1 3 4 5 7}\\
    {1 3 4 7 8}
    &{1 3 5 7 8}
    &{1 4 6 7 8}
    &{2 3 5 6 7}
    &{2 3 5 7 8}
    &{2 4 5 6 7}
    &{2 4 5 6 8}
    &{3 4 5 6 7}
    &{3 4 6 7 8}
  \end{array}
$$}

We sketch how to prove that this is, indeed, a triangulation of~$\cpp$.  As a first step we have to prove that this is
actually a combinatorial $4$-manifold.  To this end one has to verify that the link of each $k$-face is a
$(3-k)$-sphere, the crucial step being $k=0$.  Recognizing the $3$-sphere is algorithmically possible but the methods
proposed by Rubinstein, Thompson~\cite{MR95k:57015}, and others do not seem to be feasible in practice.  Fortunately,
however, for small cases, such as this example, the flip-heuristics by Bj\"orner and Lutz~\cite{MR2001h:57026} solves
the problem easily.

Since $f_2(\cpp_9)=\binom{9}{3}$ the triangulation is $3$-neighborly and hence the space is simply connected.  In view
of Theorem~\ref{inv:thm:Freedman} it then suffices to determine the intersection form.

Using elimination techniques we obtain that $$H^0(\cpp_9)\cong H^2(\cpp_9)\cong H^4(\cpp_9)\cong\ZZ$$
and all other
cohomology groups vanish.  Additionally, by keeping track of the bases transformations, one can obtain a generator for
$H^2(\cpp_9)$, for instance,
$b={1 2 5}^*
+{1 2 8}^*
+{1 3 7}^*
+{1 3 8}^*
+{1 4 6}^*
+{1 4 8}^*
-{2 3 5}^*
+{2 4 6}^*
+{2 5 6}^*
-{2 6 8}^*
-{3 4 7}^*
-{3 5 7}^*
-{3 5 8}^*
+{4 6 7}^*
-{4 7 8}^*
+{5 6 7}^*$.  Instead of computing the cup product $b\cup b$ directly, we first compute a generator for $H_2(\cpp_9)$ as
$$b\cap[\cpp_9]={5 6 7} - {5 6 8} + {5 7 8}  - {6 7 8}.$$
From this we obtain $Q_{\cpp_9}(b,b)=\langle
b,b\cap[\cpp_9]\rangle=1$.  In particular, $b\cup b$ generates $H^4(\cpp_9)$.  We conclude that the intersection form
$Q_{\cpp_9}$ has a $1\times1$-Gram matrix $(1)$ and that $$H^*(\cpp_9)\cong\ZZ[b]/(b^3).$$

The homology $2$-cycle $b\cap[\cpp_9]$ corresponds to the boundary of the non-face $3$-simplex ${5 6 7 8}$, that is,
$H_2(\cpp_9)$ is generated by an embedded $2$-sphere.  Note that this is typical in the sense that every element of the
second homology group of a closed, oriented, smooth $4$-manifold is generated by an embedded surface;
see~\cite[Proposition~1.2.3 and Chapter~2]{MR2000h:57038}.  In our case, $b\cap[\cpp_9]$ corresponds to a complex
projective line.

Now we turn to the computation of the Stiefel-Whitney homology classes of $\cpp_9$.  By \cite[Proposition
1.4.18]{MR2000h:57038} the result $Q_{\cpp_9}(b,b)\ne0$ above already implies that $w_2$ does not vanish.  A direct
computation can be sketched as follows.  Enumerating all regular pairs yields that $\omega_1$ and $\omega_3$ vanish,
$\omega_0={0} + {5} + {8}$, and that $\omega_2={0 1 2}+{0 1 8}+{0 2 3}+{0 2 7}+{0 2 8}+{0 3 4}+{0 4 5}+{0 5 6}+{0 5
  7}+{0 5 8}+{0 6 7}+{0 7 8}+{1 2 6}+{1 2 7}+{1 2 8}+{1 4 6}+{1 4 7}+{2 3 8}+{2 4 6}+{2 4 7}+{2 5 7}+{2 5 8}+{3 4 8}+{3
  5 7}+{3 5 8}+{3 7 8}+{4 5 6}+{4 5 7}+{4 5 8}+{4 6 7}+{5 6 7}+{5 6 8}+{5 7 8}+{6 7 8}$.  Reducing modulo boundaries
simplifies the result to $\omega_2={5 6 7}+{5 6 8}+{5 7 8}+{6 7 8}$ (and $\omega_0$ is the $\ZZ/2$-homology class of a
point).

\section{Fundamental Group}

For an arbitrary topological space $\cT$ and a point $*\in\cT$ let $\pi_1(\cT,*)$ be the set of homotopy classes of
closed paths based at~$*$.  The concatenation of paths imposes a group multiplication on $\pi_1(\cT,*)$.  This gives the
\emph{fundamental group} of~$\cT$ with base point~$*$.  For pathwise connected $\cT$ the structure of the fundamental
group does not depend on the choice of the base point.  In this case we omit the base point in the notation.  The space
$\cT$ is \emph{simply connected} if it is pathwise connected and $\pi_1(\cT)$ is the trivial group.

Now let $\Delta$ be a connected finite simplicial complex, and let $S\subset\Delta^{(1)}$ be the set of edges of a fixed
spanning tree. We define a finitely presented group $E(\Delta,S)$ by listing generators and relations: Each edge
$v_0v_1\in\Delta^{(1)}\setminus S$ is a generator, and for each $2$-face $w_0w_1w_2\in\Delta^{(2)}$ we have a relation
$$(w_0w_1) (w_1w_2) (w_0w_2)^{-1}=1,$$
with the understanding that for $w_iw_k\in S$ the corresponding factor is the
identity.  Then $E(\Delta,S)\cong\pi(\Delta)$.  For a combinatorial approach to fundamental groups see Seifert and
Threlfall~\cite[Chapter~7]{MR82b:55001}.

It is easy to show that each finitely presented group arises as the fundamental group of a finite simplicial complex.
In fact, it suffices to suitably modify the construction in Section~\ref{inv:homology:examples} of a simplicial complex
with a given finitely generated abelian group as its first homology; see Stillwell~\cite[\S3.4.4]{MR94a:57001}. Each
finitely presented group even arises as the fundamental group of a $4$-manifold.

Therefore, in spite of the fact that it is simple to obtain a finite description for the fundamental group of~$\Delta$,
this is often not very useful: In general, it is not possible to decide whether a given finitely presented group is
finite or trivial.  Moreover, it is, in general, impossible to decide whether two finitely presented groups are
isomorphic.  As a consequence it is impossible to decide whether two finite simplicial complexes are homotopy equivalent
or homeomorphic.  For these insolvability results see~\cite[Chapter~9]{MR94a:57001}.

However, there is an algorithm, due to Todd and Coxeter~\cite{ToddCoxeter}, which, for a given finitely presented group
$G$ and a subgroup $U$ of finite index, enumerates all cosets $Ug$ with $g\in G$.  In particular, if $G$ happens to be
finite, then the Todd-Coxeter algorithm can be used to determine the order of~$G$.  This and related algorithms can be
used to obtain simplified presentations in particular cases.

The existence of coverings is tightly related to the structure of the fundamental group of the base space;
see~\cite[Chapter~III]{MR2000b:55001} and, for a constructive approach, Rees and Soicher~\cite{MR2001h:57004}.

The first homology group $H_1(\Delta)$ coincides with the quotient of $\pi_1(\Delta)$ by its commutator group.  In
particular, if $\pi_1(\Delta)$ is abelian then $\pi_1(\Delta)\cong H_1(\Delta)$.  In this case the Todd-Coxeter
algorithm can, of course, be replaced by any Smith-Normal-Form algorithm.

Brown~\cite{MR18:753a} gave an algorithm to compute the higher homotopy groups of any simply connected finite simplicial
complex; see also Sergeraert~\cite[Theorem~10.2]{MR95c:55017}.

\section{Data Types and Other Combinatorial Models for Spaces}\label{inv:sec:data}

So far we chose a finite simplicial complex as our model of a topological space.  This allowed to access several
algebraic invariants algorithmically.  Two questions arise:
\begin{enumerate}
\item Is there a more efficient description?
\item Is there a computationally feasible model which comprises more spaces?
\end{enumerate}

The first question is easily motivated by looking at the boundary $\partial\Delta_d$ of a single $d$-simplex: Since
$f_k(\partial\Delta_d)=\binom{d+1}{k+1}$ the size of the boundary matrices is not bounded by a polynomial in~$d$, and
hence the time required to compute the simplicial homology (via the methods sketched in this survey) is exponential
in~$d$.  On the other hand, $\partial\Delta_d$ is homeomorphic to $\Sph^{d-1}$, which has a cell decomposition with two
cells only.  Applying the same Smith-Normal-Form algorithms as for the simplicial theory to the cellular boundary
matrices allows to compute all cellular homology modules of $\partial\Delta_d$ in linear time.  For an introduction to
cellular homology see~\cite[\S{}IV.10]{MR2000b:55001}.

The problem here is to determine the boundary matrices.  Since the characteristic maps, which describe how (the boundary
of) a particular $k$-cell is glued to the $(k-1)$-skeleton, can be wild, there is no hope for a purely combinatorial
description of arbitrary cell complexes.  So one solution could be to represent a cell complex by its cellular boundary
matrices.  However, this means to loose (quite a lot of) information: In general, it is not even possible to define the
cup product from the cellular boundary matrices.

There are special types of cell complexes with restricted characteristic maps, which have a useful combinatorial
description: A regular cell decomposition of a surface can be represented, e.g., as a quad-edge data structure.  For an
overview see O'Rourke~\cite[\S{}4.4.3]{MR99k:68198}.  For a certain generalization to higher dimensions see
Brisson~\cite{MR93m:68163}.

One other useful special class of regular cell complexes are the pseudo-simplicial complexes: each cell is a simplex,
but any two simplices are allowed to share several faces in their boundary.  For instance, a sphere of arbitrary
dimension~$d$ can be represented by two $d$-simplices which share their whole boundary.

Since every cell complex can be triangulated, replacing finite simplicial complexes by general finite cell complexes
does not enlarge the class of spaces accessible.  Employing techniques from functional programming Sergeraert
suggested~\cite{MR95c:55017} to encode a space~$\cT$ as a function which, for each~$k$, computes the $k$-th homology
group.  His class of \emph{spaces with effective homology}, e.g., is closed under forming loop spaces.  This way it is
feasible to compute the higher homotopy groups of simply connected finite simplicial
complexes~\cite[Theorem~10.2]{MR95c:55017}; see also~\cite{MR18:753a}.

\section{Software}

A program to compute simplicial homology is a (proposed) GAP~\cite{GAP4} package due to Dumas et
al.~\cite{gap-homology}.  As a special feature this program offers a variety of algorithms, including modular ones.

De Silva's PLEX~\cite{plex} is a collection of MATLAB functions to be used as a research tool for building and studying
simplicial complexes.  This includes $\ZZ/2$-homology computation as well as rational homology (with floating point
arithmetic).

SnapPea is a program for creating and studying hyperbolic $3$-manifolds~\cite{SnapPea}.  It focuses on metric properties.

The flip-heuristics of Bj\"orner and Lutz~\cite{MR2001h:57026} to simplify a given triangulation of a space have been
implemented by Lutz~\cite{bistellar}.

\texttt{polymake} has originally been designed as a software package for studying convex
polytopes~\cite{polymake,SoCG01:polymake}.  Yet, starting with version 2.0, \texttt{polymake}'s application
\texttt{TOPAZ} deals with finite simplicial complexes.  Besides many standard and a few non-standard constructions, it
primarily supports to compute algebraic invariants: simplicial homology and cohomology modules (with and without bases,
via elimination), cap-products, cup-products, Stiefel-Whitney classes, and intersection forms of $4$-manifolds.  In
fact, the example computation for the complex projective plane, see Section~\ref{inv:subsec:cpp}, had been performed
with \texttt{polymake/TOPAZ}.


\bibliographystyle{plain}
\bibliography{invariants}

\end{document}